\documentstyle{article}

\def\esp{\raisebox{.3ex}{!`}}

\def\pr{\mbox{\small \O}}

\def\Rel{{\it Rel}}
\def\Fun{{\it Fun}}

\def\Mstr{${\cal M}^{str}$}
\def\LL{${\cal L}^L$}
\def\LR{${\cal L}^R$}
\def\Mc{${\cal M}_c$}
\def\Lc{${\cal L}_c$}
\def\Mcstr{${\cal M}_c^{str}$}
\def\Lcstr{${\cal L}_c^{str}$}
\def\Rmin{${\cal R}^-$}
\def\Rstr{${\cal R}^{str}$}
\def\Cstr{${\cal C}^{str}$}

\def\mj{{\mathbf{1}}}

\def\cirk{\,{\raisebox{.3ex}{\tiny $\circ$}}\,}
\def\prop#1#2{\vspace{2ex} \noindent{\sc #1.} {\it #2} \par \vspace{2ex}}
\def\dkz{\noindent{\sc Proof. }}
\def\qed{\hfill $\dashv$}
\def\str{\rightarrow}

\begin{document}

\title{Coherence for Monoidal Endo\-functors}
\author{\small {\sc Kosta Do\v sen} and {\sc Zoran Petri\' c}
\\[1ex]
{\small Mathematical Institute, SANU}\\[-.5ex]
{\small Knez Mihailova 36, p.f.\ 367, 11001 Belgrade,
Serbia}\\[-.5ex]
{\small email: \{kosta, zpetric\}@mi.sanu.ac.rs}}
\date{}
\maketitle

\begin{abstract}
\noindent The goal of this paper is to prove coherence results
with respect to relational graphs for monoidal endo\-functors,
i.e.\ endo\-functors of a monoidal category that preserve the
monoidal structure up to a natural transformation that need not be
an isomorphism. These results are proved first in the absence of
symmetry in the monoidal structure, and then with this symmetry.
In the later parts of the paper the coherence results are extended
to monoidal endo\-functors in monoidal categories that have
diagonal or codiagonal natural transformations, or where the
monoidal structure is given by finite products or coproducts.
Monoidal endo\-functors are interesting because they stand behind
monoidal monads and comonads, for which coherence will be proved
in a sequel to this paper.
\end{abstract}

\noindent {\small \emph{Mathematics Subject Classification
(2000):} 18D10, 18C05, 18A15, 03F07, 03F05}

\vspace{.5ex}

\noindent {\small {\it Keywords:} monoidal endo\-functor,
coherence, relational graphs, finite products, finite coproducts}

\section{Introduction}
A monoidal functor is a functor between monoidal categories that
preserves the monoidal structure up to a natural transformation
that need not be an isomorphism (see Section~3 below; this notion
stems from \cite{EK66}, Section II.1). Coherence results for
monoidal functors were obtained long ago in \cite{Ep66} and
\cite{L72}. In \cite{Ep66} one can find a result for such functors
between symmetric monoidal categories in the absence of unit
objects, while in \cite{L72} unit objects are allowed, and
nonsymmetric monoidal categories are considered too.

To get coherence with the unit objects, \cite{L72} introduces
implicitly graphs that connect occurrences of the generating
functor (see the beginning of the next section below). The
standard graphs, which stem from \cite{KML71}, and earlier work of
Mac Lane and Kelly, connect occurrences of generating objects.

Our goal in this paper is first to extend these old coherence
results to the situation where we have not a monoidal functor
between two categories, but an endo\-functor of a single monoidal
category. This involves matters that go beyond \cite{L72}, where
application of functors cannot be iterated. Monoidal
endo\-functors are interesting because they stand behind
\emph{monoidal} monads and comonads, and the present paper lays
the ground for a study of coherence in these monads and comonads.

Monoidal monads stem from \cite{Ko70} and \cite{Ko72}. More recent
papers on monoidal comonads are \cite{Moe02}, \cite{BV07} and
\cite{PS09}. We will prove coherence results for monoidal monads
and comonads in a sequel to this paper \cite{DP09b}. In the
present paper, and in that sequel, we understand coherence with
respect to graphs that are like those of \cite{L72}. Coherence
states that there is a faithful functor from a freely generated
categorial structure, for which we prove coherence, into the
category whose arrows are such graphs.  We obtain thereby a
characterization of the freely generated categorial structure in
terms of graphs. Such coherence results give very useful
procedures for deciding whether a diagram of canonical arrows
commutes. (A general treatment of coherence in this spirit may be
found in~\cite{DP04}.)

One finds in \cite{Mog91} a notion of monad inspired by
\cite{Ko70} and \cite{Ko72}, at the basis of which one finds the
notions of left and right monoidal endo\-functors, for which we
are also going to prove coherence. We prove our results first in
the absence of symmetry, and then with symmetry. In the later part
of the paper we extend our coherence results to monoidal
endo\-functors in monoidal categories that have diagonal or
codiagonal natural transformations (we call these monoidal
categories \emph{relevant} categories), or where the monoidal
structure is given by finite products or coproducts.

Some of our coherence results may be understood as basic coherence
results for equations between deductions in modal logic. In this
paper we find systems that may be understood as fragments of $K$
with the necessity operator $\Box$ primitive; in the sequel, with
comonads, we will find fragments of $S4$ with $\Box$ primitive.

\section{Endo\-functors in monoidal categories}

In this section we deal with coherence for monoidal categories
with endo\-functors for which we do not assume yet that they are
monoidal. This is a basic auxiliary result, which we will need
later.

A monoidal category is, as usual, a category with a
biendo\-functor $\otimes$, a special unit object $I$, and the
natural isomorphisms whose components are the arrows

\begin{tabbing}
\hspace{9em}\= $a_{A,B,C}\!:(A\otimes B)\otimes C\str
A\otimes(B\otimes C)$,
\\[1.5ex]
\> $l_A\!:I\otimes A\str A$,\hspace{3em}$r_A\!:A\otimes I\str A$,
\end{tabbing}
which satisfy Mac Lane's coherence equations (see \cite{ML98},
Section VII.1; our notation comes from \cite{EK66}, Section II.1).

Let $\cal E$ be the free monoidal category with a family of
endo\-functors; freedom means here and later free generation by
two arbitrary sets, one of which is conceived as the set of
\emph{generating objects}, and the other as the set of
\emph{generating functors} $\{E^i\mid i\in {\cal I}\}$. We call
the generating objects and the generating functors collectively
\emph{generators}. If $\cal I$ is empty, then $\cal E$ is just the
free monoidal category generated by a set of generating objects.

The category $\cal E$ is made of syntactical material. Its objects
are \emph{propositional formulae} built with the binary connective
$\otimes$, the unary connectives $E^i$ and the nullary connective
$I$ out of the generating objects, which we take to be the
propositional letters $p,q,r,\ldots$ An object of $\cal E$ is
\emph{atomic} when it is a generating object or of the form
$E^iA$. An object of $\cal E$ is \emph{diversified on generating
objects} when every generating object occurs in it at most once.
We define analogously diversification on generating functors, and
we say that an object is \emph{diversified} when it is diversified
both on generating objects and on generating functors. For $E^iA$
a subformula of an object $B$ of $\cal E$, the \emph{scope in} $B$
of the outermost occurrence of $E^i$ in $E^iA$ is the set of all
the generators in $A$.

The arrows of $\cal E$ are equivalence classes of arrow terms made
out of the primitive arrow terms $\mj_A$, $a_{A,B,C}$, $l_A$,
$r_A$, with the operations $\cirk$, $\otimes$ and $E^i$ so that
the equations assumed for defining a monoidal category together
with functorial equations for $E^i$, for each $i\in\cal I$, are
satisfied (cf.\ \cite{DP04}, Chapter~2). We take for granted the
superscripts $i$ of $E^i$ and omit them, except when they are
essential. (We do the same later with $\psi$, $\psi_0$, $\psi^L$
and $\psi^R$.) The existence of free structures like $\cal E$ is
guaranteed by the purely equational definition of these
structures.

Every arrow term $f$ of $\cal E$ is equal to an arrow term
$f_n\cirk\ldots\cirk f_1$, called \emph{developed}, which is
$\mj_A$ if $n=0$, and if $n\geq 1$, then for every
$j\in\{1,\ldots,n\}$ in $f_j$ we have exactly one occurrence of
$a$, $l$ or $r$, and no occurrence of $\cirk$; such an $f_i$ is
called a \emph{factor}. The subterm $a_{A,B,C}$, or $l_A$, or
$r_A$, of a factor is its \emph{head} (see \cite{DP04}, Section
2.7). A factor with head $a_{A,B,C}$ is an $a$-\emph{factor}, and
analogously in other cases. We are going to prove the following
theorem.

\prop{$\cal E$-Coherence}{The category $\cal E$ is a preorder.}

\dkz Suppose we have two arrow terms $f,g\!:A\str B$ of $\cal E$.
To show that $f=g$, we proceed by induction on the number $n$ of
occurrences of $E$ in $A$, which is equal to this number in $B$.
In the basis, when $n=0$, we have Mac Lane's coherence result for
monoidal categories (see \cite{ML98}, Section VII.2, or
\cite{DP04}, Chapter~4).

When $n>0$, take a single arbitrary occurrence of $E^i$ in $A$,
and replace it by $E^j$ such that $j$ is not an index of any $E$
in $A$. (If all the generating functors occur in $A$, then we
enlarge for the sake of the proof the set of generating functors
with a new functor $E^j$, which functions just as a placeholder.)
Make this replacement at the appropriate place in $B$, and in the
arrow terms $f$ and $g$ so as to obtain the arrow terms
$f',g'\!:A'\str B'$ of $\cal E$. By naturality and functorial
equations $f'$ is equal to a developed arrow term $f_2\cirk f_1$
such that no head of a factor of $f_1\!:A'\str C$ is in the scope
of $E^j$ and all the heads of the factors of $f_2\!:C\str B'$ are
in the scope of $E^j$. The object $C$ is completely determined by
$A'$ and $B'$. Analogously we have $g'=g_2\cirk g_1$, with
$g_1\!:A'\str C$ and $g_2\!:C\str B'$. Since $E^jD$ in $f_1$ and
$g_1$ amounts to a generating object (it is a parameter), and
since in $f_2$ and $g_2$ only what is within the scope of $E^j$
counts, by the induction hypothesis we have $f_1=g_1$ and
$f_2=g_2$, and so $f'=g'$, from which $f=g$ follows by
substitution. \mbox{\hspace{1em}}\qed

\section{Monoidal and locally monoidal endo\-functors}

An endo\-functor may preserve the monoidal structure of a category
up to a natural transformation either globally or locally. We have
\emph{global} preservation when our endo\-functor $E$ is a
\emph{monoidal functor} in the sense of \cite{EK66} (Section II.1;
see also \cite{ML98}, Section XI.2). This means that in our
monoidal category we have a natural transformation whose
components are the arrows
\[
\psi_{A,B}\!:EA\otimes EB\str E(A\otimes B),
\]
and we have also the arrow $\psi_0\!:I\str EI$; the monoidal
structure is preserved up to $\psi$ and $\psi_0$, which means that
the following equations hold:
\begin{tabbing}
\hspace{1.5em}\=${(\psi a)}$\hspace{3em}\=$Ea_{A,B,C}\cirk
\psi_{A\otimes B,C}\cirk(\psi_{A,B}\otimes\mj_{EC})=$\kill

${(\psi a)}$\hspace{.8em}$Ea_{A,B,C}\cirk \psi_{A\otimes
B,C}\cirk(\psi_{A,B}\otimes\mj_{EC})=\psi_{A,B\otimes
C}\cirk(\mj_{EA}\otimes\psi_{B,C})\cirk a_{EA,EB,EC}$,
\\[2ex]
\>${(\psi
l)}$\>$El_A\cirk\psi_{I,A}\cirk(\psi_0\otimes\mj_{EA})\;$\=$=l_{EA}$,
\\[2ex]
\>${(\psi
r)}$\>$Er_A\cirk\psi_{A,I}\cirk(\mj_{EA}\otimes\psi_0)$\>$=r_{EA}$.
\end{tabbing}
The global character of the preservation is manifested in ${(\psi
a)}$ by $E$ from the left-hand side falling on every index of $a$
on the right-hand side. (The notation with $\psi$ stems from
\cite{Ko70} and \cite{Ko72}.)

We have \emph{local} preservation with the following three notions
of endo\-functor suggested by \cite{Ko70} and \cite{Ko72}.
Monoidal functors need not be endo\-functors, but the notions we
are going to consider now are tied to endo\-functors only.

We say that an endo\-functor $E$ of a monoidal category is
\emph{left monoidal} when we have a natural transformation whose
components are the arrows
\[
\psi^L_{A,B}\!:EA\otimes B\str E(A\otimes B),
\]
and the monoidal structure is preserved up to $\psi^L$, which
means that the following equations hold:
\begin{tabbing}
\hspace{1.5em}\=${(\psi^L a)}$\hspace{3em}\=$Ea_{A,B,C}\cirk
\psi^L_{A\otimes
B,C}\cirk(\psi^L_{A,B}\otimes\mj_C)=\psi^L_{A,B\otimes C}\cirk
a_{EA,B,C}$,
\\[2ex]
\>${(\psi^L r)}$\>$Er_A\cirk\psi^L_{A,I}=r_{EA}$.
\end{tabbing}
The local character of the preservation is manifested in ${(\psi^L
a)}$ by $E$ from the left-hand side falling on a single index of
$a$ on the right-hand side.

We say, analogously, that $E$ is \emph{right monoidal} when we
have a natural transformation whose components are the arrows
\[
\psi^R_{A,B}\!:A\otimes EB\str E(A\otimes B),
\]
and the monoidal structure is preserved up to $\psi^R$, which
means that the following equations hold:
\begin{tabbing}
\hspace{1.5em}\=${(\psi^R a)}$\hspace{3em}\=$Ea_{A,B,C}\cirk
\psi^R_{A\otimes B,C}=\psi^R_{A,B\otimes
C}\cirk(\mj_A\otimes\psi^R_{B,C})\cirk a_{A,B,EC}$,
\\[2ex]
\>${(\psi^R l)}$\>$El_A\cirk\psi^R_{I,A}=l_{EA}$.
\end{tabbing}

We say that $E$ is \emph{locally} monoidal when it is both left
and right monoidal, and we have moreover the equation
\begin{tabbing}
${(\psi^L\psi^R a)}$\hspace{1.2em}$Ea_{A,B,C}\cirk\psi^L_{A\otimes
B,C}\cirk(\psi^R_{A,B}\otimes\mj_C)=\psi^R_{A,B\otimes
C}\cirk(\mj_A\otimes\psi^L_{B,C})\cirk a_{A,EB,C}$.
\end{tabbing}

We define, as we defined the category $\cal E$ in the preceding
section, the free monoidal categories with a family of monoidal
endo\-functors, a family of left monoidal endo\-functors, a family
of right monoidal endo\-functors, or a family of locally monoidal
endo\-functors, which we call respectively $\cal M$, \LL, \LR\ and
$\cal L$. All these categories have the same propositional
formulae as objects (provided the sets of generators are the
same). Then it is easy to see that the category \LR\ is an
isomorphic, mirror image, of \LL. For the categories $\cal M$,
\LL, \LR\ and $\cal L$, we define the notions of developed arrow
term, factor and head of a factor analogously to what we had for
$\cal E$ in the preceding section.

We define a functor $G$ from $\cal M$ to the category \Fun\ of
functions between finite ordinals by stipulating that $GA$, for
$A$ a propositional formula, is the number of occurrences of $E$
in $A$ (i.e.\ the number of all $E^i$'s in $A$, for every $i$),
while $Gf$ for an arrow term $f$ of $\cal M$ is defined
inductively on the complexity of $f$. We have that $G\mj_A$,
$Ga_{A,B,C}$, $Gl_A$ and $Gr_A$ are identity functions, while for
the remaining primitive arrow terms we have the clause
corresponding to the following picture:
\begin{center}
\begin{picture}(30,50)
\put(2,16){\line(0,1){20}}

\put(12,16){\line(0,1){20}} \put(16,16){\line(0,1){20}}

\put(11.75,16){\line(1,0){4.3}} \put(11.75,36){\line(1,0){4.3}}

\put(40,16){\line(0,1){20}} \put(44,16){\line(0,1){20}}

\put(39.75,16){\line(1,0){4.3}} \put(39.75,36){\line(1,0){4.3}}

\put(2,16){\line(3,2){30}}

\put(0,42){\makebox(0,0)[l]{$E\;A\otimes E\,B$}}
\put(-1.75,8){\makebox(0,0)[l]{$E\,(A\;\otimes\;\,B)$}}

\put(-45,25){\makebox(0,0)[l]{$G\psi_{A,B}$}}
\end{picture}
\end{center}
and $G\psi_0$ is the empty function from $\pr$, which is $GI$, to
$\{\pr\}$, which is $GEI$ (in our picture we obtain a crossing if
there is an $E$ in $A$, different or not from the $E$ in the
picture). We also have clauses corresponding to the following
pictures:
\begin{center}
\begin{picture}(120,30)(0,10)

\put(10,16){\line(0,1){20}} \put(30,16){\line(0,1){20}}

\put(10,16){\line(1,0){20}} \put(10,36){\line(1,0){20}}

\put(40,16){\line(0,1){20}} \put(60,16){\line(0,1){20}}

\put(40,16){\line(1,0){20}} \put(40,36){\line(1,0){20}}

\put(14,25){\makebox(0,0)[l]{$Gf$}}
\put(44,25){\makebox(0,0)[l]{$Gg$}}

\put(-45,25){\makebox(0,0)[l]{$G(f\otimes g)$}}

\put(150,16){\line(0,1){20}} \put(170,16){\line(0,1){20}}

\put(150,16){\line(1,0){20}} \put(150,36){\line(1,0){20}}

\put(140,16){\line(0,1){20}}

\put(154,25){\makebox(0,0)[l]{$Gf$}}

\put(105,25){\makebox(0,0)[l]{$GEf$}}

\end{picture}
\end{center}
and  $G(g\cirk f)$ is the composition of the functions $Gf$ and
$Gg$. It is easy to verify by induction on the length of
derivation that $G$ so defined on the arrow terms of $\cal M$
induces a functor from $\cal M$ to \Fun.

Intuitively, with $Gf\!:GA\str GB$ we note from which occurrences
of $E$ in $A$ the occurrences of $E$ in $B$ originate. We call
$Gf$ a \emph{graph}. For example, for the two sides of the
equation ${(\psi a)}$ we have the following pictures:
\begin{center}
\begin{picture}(230,110)
\put(3,17){\line(0,1){18}} \put(3,47){\line(0,1){18}}
\put(3,77){\line(0,1){18}} \put(63,77){\line(0,1){18}}

\put(3,47){\line(3,1){55}} \put(3,77){\line(3,2){27}}

\put(-3,100){\makebox(0,0)[l]{$(EA\otimes EB)\otimes EC$}}
\put(0,70){\makebox(0,0)[l]{$E(A\otimes B)\,\otimes\, EC$}}
\put(0,40){\makebox(0,0)[l]{$E((A\otimes B)\otimes C)$}}
\put(0,10){\makebox(0,0)[l]{$E(A\otimes(B\otimes C))$}}

\put(203,17){\line(0,1){18}} \put(203,47){\line(0,1){18}}
\put(235,47){\line(0,1){18}} \put(203,77){\line(0,1){18}}
\put(235,77){\line(0,1){18}} \put(267,77){\line(0,1){18}}

\put(235,47){\line(3,2){27}} \put(203,17){\line(3,2){27}}

\put(197,100){\makebox(0,0)[l]{$(EA\,\otimes\, EB)\otimes EC$}}
\put(200,70){\makebox(0,0)[l]{$EA\otimes (EB\,\otimes\, EC)$}}
\put(200,40){\makebox(0,0)[l]{$EA\,\otimes\, E(B\otimes C)$}}
\put(200,10){\makebox(0,0)[l]{$E(A\otimes(B\otimes C))$}}

\put(-10,85){\makebox(0,0)[r]{$\psi_{A,B}\otimes\mj_C$}}
\put(-10,55){\makebox(0,0)[r]{$\psi_{A\otimes B,C}$}}
\put(-10,25){\makebox(0,0)[r]{$Ea_{A,B,C}$}}

\put(190,55){\makebox(0,0)[r]{$\mj_{EA}\otimes\psi_{B,C}$}}
\put(190,25){\makebox(0,0)[r]{$\psi_{A,B\otimes C}$}}
\put(190,85){\makebox(0,0)[r]{$a_{EA,EB,EC}$}}

\end{picture}
\end{center}
and for the two sides of ${(\psi l)}$ we have the following
pictures:

\begin{center}
\begin{picture}(180,110)
\put(3,17){\line(0,1){18}} \put(3,47){\line(0,1){18}}
\put(33,77){\line(0,1){18}}

\put(3,47){\line(3,2){27}}

\put(12,100){\makebox(0,0)[l]{$I\otimes EA$}}
\put(0,70){\makebox(0,0)[l]{$EI\,\otimes\, EA$}}
\put(0,40){\makebox(0,0)[l]{$E(I\otimes A)$}}
\put(0,10){\makebox(0,0)[l]{$EA$}}

\put(153,77){\line(1,1){18}}

\put(150,100){\makebox(0,0)[l]{$I\otimes EA$}}
\put(148,70){\makebox(0,0)[l]{$EA$}}

\put(-10,85){\makebox(0,0)[r]{$\psi_0\otimes\mj_{EA}$}}
\put(-10,55){\makebox(0,0)[r]{$\psi_{I,A}$}}
\put(-10,25){\makebox(0,0)[r]{$El_A$}}

\put(140,85){\makebox(0,0)[r]{$l_{EA}$}}

\end{picture}
\end{center}

The functors from \LL, \LR\ and $\cal L$ to \Fun\ analogous to
$G$, which we all call $G$, are defined as $G$ save for the
clauses corresponding to the following pictures:
\begin{center}
\begin{picture}(160,50)
\put(2,16){\line(0,1){20}}

\put(12,16){\line(0,1){20}} \put(16,16){\line(0,1){20}}

\put(11.75,16){\line(1,0){4.3}} \put(11.75,36){\line(1,0){4.3}}

\put(40,16){\line(0,1){20}} \put(44,16){\line(0,1){20}}

\put(39.75,16){\line(1,0){4.3}} \put(39.75,36){\line(1,0){4.3}}

\put(0,42){\makebox(0,0)[l]{$E\;A\;\,\otimes \;\,B$}}
\put(-1.75,8){\makebox(0,0)[l]{$E\,(A\;\otimes\;\,B)$}}

\put(-45,25){\makebox(0,0)[l]{$G\psi^L_{A,B}$}}

\put(142,16){\line(0,1){20}} \put(146,16){\line(0,1){20}}

\put(141.75,16){\line(1,0){4.3}} \put(141.75,36){\line(1,0){4.3}}

\put(170,16){\line(0,1){20}} \put(174,16){\line(0,1){20}}

\put(169.75,16){\line(1,0){4.3}} \put(169.75,36){\line(1,0){4.3}}

\put(132,16){\line(3,2){30}}

\put(130,42){\makebox(0,0)[l]{$\;\;\;\;A\otimes E\,B$}}
\put(128.25,8){\makebox(0,0)[l]{$E\,(A\;\otimes\;\,B)$}}

\put(85,25){\makebox(0,0)[l]{$G\psi^R_{A,B}$}}

\end{picture}
\end{center}
(this means that $G\psi^L_{A,B}$ is an identity arrow). The target
category of $G$ for \LR\ and $\cal L$ is the subcategory of \Fun\
of bijections between finite ordinals, and for \LL\ this is the
discrete subcategory of \Fun, with all arrows just identity
arrows.

If $\cal K$ is a category like $\cal M$, \LL, \LR\ or $\cal L$,
then we call $\cal K$-\emph{Coherence} the proposition that $G$
from $\cal K$ to \Fun, or a category like \Fun, is a faithful
functor. Since the image of \LL\ under $G$ is a discrete category,
\LL-Coherence amounts to the proposition that \LL\ is a preorder,
and since \LL\ and \LR\ are isomorphic, \LR-Coherence amounts too
to the proposition that \LR\ is a preorder. (Our notion of $\cal
K$-coherence is a standard notion of coherence, which stems from
Mac Lane's coherence results for monoidal and symmetric monoidal
categories; see \cite{ML98}, \cite{DP04} and references therein.)

Note that with this understanding of coherence we cannot expect
that $\psi_{A,B}$ and $\psi_0$ be isomorphisms. With the natural
$G$ image (now not in \Fun) of the inverses of $\psi_{A,B}$ and
$\psi_0$, we have
\begin{center}
\begin{picture}(180,70)(0,-5)
\put(2,36){\line(0,1){20}}

\put(12,36){\line(0,1){20}} \put(16,36){\line(0,1){20}}

\put(11.75,36){\line(1,0){4.3}} \put(11.75,56){\line(1,0){4.3}}

\put(40,36){\line(0,1){20}} \put(44,36){\line(0,1){20}}

\put(39.75,36){\line(1,0){4.3}} \put(39.75,56){\line(1,0){4.3}}

\put(2,36){\line(3,2){30}}

\put(2,3){\line(0,1){20}}

\put(12,3){\line(0,1){20}} \put(16,3){\line(0,1){20}}

\put(11.75,3){\line(1,0){4.3}} \put(11.75,23){\line(1,0){4.3}}

\put(40,3){\line(0,1){20}} \put(44,3){\line(0,1){20}}

\put(39.75,3){\line(1,0){4.3}} \put(39.75,23){\line(1,0){4.3}}

\put(2,23){\line(3,-2){30}}

\put(0,62){\makebox(0,0)[l]{$E\;A\otimes E\,B$}}
\put(-1.75,28){\makebox(0,0)[l]{$E\,(A\;\otimes\;\,B)$}}
\put(0,-4){\makebox(0,0)[l]{$E\;A\otimes E\,B$}}

\put(140,62){\makebox(0,0)[l]{$EI$}}
\put(147,28){\makebox(0,0)[l]{$I$}}
\put(140,-4){\makebox(0,0)[l]{$EI$}}

\end{picture}
\end{center}
neither of which corresponds to an identity arrow. The reasons for
this failure of isomorphism are similar to the reasons for the
failure of the isomorphism of distribution investigated in
\cite{DP04}.

Instead of formulating our coherence results in terms of $G$ and
graphs, we could have formulations based on diversified objects
(see the preceding section). For example, $\cal M$-Coherence,
which we are going to prove in the next section, is equivalent to
the proposition that for all arrow terms $f,g\!:A\str B$ of $\cal
M$ with $B$ diversified we have $f=g$ in $\cal M$.

\section{$\cal M$-Coherence}

The category $\cal M$ is equivalent to its strictification \Mstr,
where
\begin{tabbing}
\hspace{1.5em}\=$(A\otimes B)\otimes C=A\otimes(B\otimes
C)$,\hspace{5em}\=$a_{A,B,C}=\mj_{A\otimes B\otimes C}$,
\\[1.5ex]
\>$I\otimes A=A=A\otimes I$,\>$l_A=\mj_A=r_A$.
\end{tabbing}
The preordered groupoid subcategory of $\cal M$ over which we make
the strictification is the category $\cal E$ of Section~2 (see
\cite{DP04}, Section 3.2, which, together with Section 3.1,
provides a general treatment of strictification, where references
to earlier approaches may be found).

Let $H$ be the functor from \Mstr\ to $\cal M$, and $H'$ the
functor in the opposite direction, by which \Mstr\ and $\cal M$
are equivalent categories. We will show that the composite functor
$GH$ from \Mstr\ to \Fun\ is faithful. This implies that $G$ from
$\cal M$ to \Fun\ is faithful, i.e.\ $\cal M$-Coherence, in the
following manner. Suppose that $Gf=Gg$; then, since for every
arrow $h$ of $\cal M$ we have $GHH'h=Gh$, we obtain $GHH'f=GHH'g$,
and by the faithfulness of $GH$, we obtain $H'f=H'g$, from which
we obtain $HH'f=HH'g$, and hence $f=g$ in $\cal M$.

\prop{Proposition~1}{The functor $GH$ from \Mstr\ to \Fun\ is
faithful.}

\dkz Every arrow term $f$ of \Mstr\ is equal to a developed arrow
term $f_n\cirk\ldots\cirk f_1$, such that each $f_j$ is either a
$\psi$-factor or a $\psi_0$-factor. By applying naturality and
functorial equations, and the equations
\[
\psi_{I,A}\cirk(\psi_0\otimes\mj_{EA})=\mj_{EA}=\psi_{A,I}\cirk(\mj_{EA}\otimes\psi_0),
\]
which are ${(\psi l)}$ and ${(\psi r)}$ strictified, we obtain
from a developed arrow term an arrow term equal to it, of the form
$h\cirk g_m\cirk\ldots\cirk g_1$, where in $h$ we have no
occurrence of $\psi$, and $\cirk$ may occur only in subterms of
$h$ of the form $E\ldots E\psi_0\cirk\ldots\cirk\psi_0\!:I\str
E\ldots EI$; for $g_m\cirk\ldots\cirk g_1$ we assume that it is
developed without $\psi_0$-factors.

If $m\geq 1$, then each $g_j$ is a $\psi$-factor, and we assign to
$g_j$ a finite ordinal $\tau(g_j)$ obtained by applying the
function $GH(g_m\cirk\ldots\cirk g_{j+1})$ to the number
$\kappa(g_j)$ of occurrences of $E$ in $g_j$ to the left of
$\psi$. Intuitively, this is the place where the contracted $E$ of
$g_j$ will end up in the codomain of $g_m$. For example, with
$g_2\cirk g_1$ being
\[
E^1\psi^2_{p,E^1(p\otimes q)}\cirk E^1(\mj_{E^2p}\otimes
E^2\psi^1_{p,q}),
\]
$\kappa(g_1)=3$ and $\tau(g_1)=2$, which is clear from the
following picture:
\begin{center}
\begin{picture}(110,80)
\put(3,17){\line(0,1){18}} \put(3,47){\line(0,1){18}}
\put(20,17){\line(0,1){18}} \put(20,47){\line(0,1){18}}
\put(50,47){\line(0,1){18}} \put(65,47){\line(0,1){18}}

\put(20,17){\line(3,2){27}} \put(51,17){\line(2,3){11.5}}
\put(65,47){\line(3,2){27}}

\put(0,70){\makebox(0,0)[l]{$E^1(E^2p\otimes E^2(E^1p \otimes E^1
q))$}} \put(0,40){\makebox(0,0)[l]{$E^1(E^2p\otimes
E^2\,E^1(p\otimes q))$}} \put(0,10){\makebox(0,0)[l]{$E^1\,
E^2(p\otimes E^1(p\otimes q))$}}

\put(-10,55){\makebox(0,0)[r]{$GHg_1$}}
\put(-10,25){\makebox(0,0)[r]{$GHg_2$}}

\put(3,0){\makebox(0,0)[b]{\scriptsize $0$}}
\put(19,0){\makebox(0,0)[b]{\scriptsize $1$}}
\put(51,0){\makebox(0,0)[b]{\scriptsize $2$}}

\end{picture}
\end{center}

It is not difficult to see that for the $\psi$-factors $g_i$ and
$g_{i+1}$ such that $\tau(g_{i+1})=k<l=\tau(g_i)$ we have, by
naturality and functorial equations, that $g_{i+1}\cirk
g_i=g'_{i+1}\cirk g'_i$ for some $\psi$-factors $g'_i$ and
$g'_{i+1}$ such that $\tau(g'_i)=k$ and $\tau(g'_{i+1})=l$. So
$g_m\cirk\ldots\cirk g_1$ is equal to $g'_m\cirk\ldots\cirk g'_1$
such that $\tau(g'_{i+1})\geq\tau(g'_i)$. If
$\tau(g'_{i+1})=\tau(g'_i)$, then they can be permuted by the
equation ${(\psi a)}$ strictified. (Note that with five
applications of that equation we may permute also the rightmost
two factors of $\psi_{p\otimes q,r\otimes s}\cirk(\mj_{E(p\otimes
q)}\otimes\psi_{r,s})\cirk(\psi_{p,q}\otimes\mj_{Er\otimes Es})$.)
We take $h\cirk g'_m\cirk\ldots\cirk g'_1$ to be a normal form of
$f$. As an arrow term, this normal form is not unique, because we
may have differences based on the last mentioned permutations or
on equations like $E\mj_A=\mj_{EA}$.

We may show however that if $GHf=GHf'$ and $f$ and $f'$ are in
normal form, then $f$ and $f'$ differ from each other only with
respect to what is mentioned in the preceding sentence. In $f$ let
a \emph{block} $\vec{f}_i$ be a composition of $\psi$-factors
$f_{i_k}\cirk\ldots\cirk f_{i_1}$ such that
$\tau(f_{i_k})=\ldots=\tau(f_{i_1})=l$. We stipulate then that
$\tau(\vec{f}_i)=l$. Let $f$ be $h\cirk\vec{f}_n\cirk\ldots\cirk
\vec{f}_1$ such that for $i,j\in\{1,\ldots,n\}$ if $i<j$, then
$\tau(\vec{f}_i)<\tau(\vec{f}_j)$. The arrangement of these blocks
is \emph{strictly increasing}. Let analogously $f'$ be $h'\cirk
\vec{f}'_n\cirk\ldots\cirk \vec{f}'_1$. From $GHf=GHf'$ we
conclude first that $n=n'$, and we proceed by induction on $n$. If
$n=0$, then we conclude easily that $h=h'$. If $n>0$, then we
conclude that $\tau(\vec{f}_1)=\tau(\vec{f}'_1)$, that
$GH\vec{f}_1=GH\vec{f}'_1$ and that
$GH(h\cirk\vec{f}_n\cirk\ldots\cirk \vec{f}_2)=GH(h'\cirk
\vec{f}'_n\cirk\ldots\cirk \vec{f}'_2)$. From this we conclude
that $\vec{f}_1=\vec{f}'_1$, and by the induction hypothesis
$h\cirk\vec{f}_n\cirk\ldots\cirk \vec{f}_2=h'\cirk
\vec{f}'_n\cirk\ldots\cirk \vec{f}'_2$. So $f=g$ in \Mstr. \qed

\vspace{2ex}

\noindent Hence we have $\cal M$-Coherence.

Note that $\cal M$ is not a preorder. The following two arrows:
\[
(\psi_0\otimes\mj_{EI})\cirk
l^{-1}_{EI},(\mj_{EI}\otimes\psi_0)\cirk r^{-1}_{EI}\!:EI\str
EI\otimes EI,
\]
have different $G$ images, and are different in $\cal M$; this
counterexample for preorder is from \cite{L72} (Section 0).
Another counterexample is given by the two arrows
\[
(\mj_{EA}\otimes(El_B\!\cirk\psi_{I,B}))\!\cirk
a_{EA,EI,EB},(Er_A\!\cirk\psi_{A,I})\otimes\mj_{EB}\!:(EA\otimes
EI)\otimes EB\!\str\! EA\otimes EB.
\]
This counterexample shows that graphs are essential for coherence
even in the absence of $\psi_0$.

Let, however, ${\cal M}^-$ be the category defined like $\cal M$
save that we reject $I$ and everything that involves it---namely,
$l$, $r$ and $\psi_0$. The category ${\cal M}^-$ is a preorder,
and graphs are irrelevant for its coherence. When we try to
determine whether there is an arrow of ${\cal M}^-$ of a given
type (i.e.\ with a given source and target), we find that if there
is such an arrow it must be unique.

This will become clear with the following example. Suppose we want
to determine whether there is an arrow
\[
f\!:E(Ep\otimes E(Eq\otimes Ep))\str EE(p\otimes E(q\otimes p)).
\]
We diversify first the propositional letters and the occurrences
of $E$ in the target, and the question is then whether we have an
arrow
\[
f'\!:E(Ep\otimes E(Eq\otimes Er))\str E^1E^2(p\otimes E^3(q\otimes
r))
\]
for some superscripts assigned to the occurrences of $E$ in the
source. The leftmost $E$ in the source must be $E^1$. Since this
$E$ has $\{p,q,r\}$ in its scope as $E^1$ in the target, we are
done with $E^1$. The leftmost of the remaining $E$'s in the source
must be $E^2$. Since this $E$ has only $\{p\}$ in its scope, while
$E^2$ in the target has $\{p,q,r\}$, we take as $E^2$ the leftmost
of the remaining $E$'s in the source in whose scope we find
$\{q,r\}$. By iterating this procedure we find the arrow
\[
E^1E^2(\mj_p\otimes\psi^3_{a,r})\cirk E^1\psi^2_{p,E^3q\otimes
E^3r}\!:E^1(E^2p\otimes E^2(E^3q\otimes E^3r))\str E^1E^2(p\otimes
E^3(q\otimes r)).
\]

\section{\LL, \LR\ and $\cal L$-Coherence}

To prove \LL-Coherence, which as we said towards the end of
Section~3 amounts to \LL\ being a preorder, we proceed as for
$\cal M$-Coherence. We introduce the strictification of \LL\ and
in the proof of the faithfulness of $GH$ we have a normal form
that is a simplified version of the normal form of the preceding
proof for \Mstr. The $h$ part of the normal form with
$\psi_0$-factors does not exist, and instead of the $\psi$-factors
part we have a $\psi^L$-factors part. The blocks are now of length
1, i.e.\ single $\psi^L$-factors, according to the equation
${(\psi^L a)}$ strictified, and $\psi^L$-factors with
$\psi^L_{A,I}$ are identity arrows by the equation ${(\psi^L r)}$
strictified. We may prove \LR-Coherence either directly in the
same manner, or just appeal to the isomorphism of \LL\ and \LR.

To prove $\cal L$-Coherence we proceed again as before. The $h$
part of the normal form does not exist again, and we have
$\psi^L$-factors and $\psi^R$-factors. A block is either a single
$\psi^L$-factor, or a single $\psi^R$-factor, or a pair of factors
made of one $\psi^L$-factor and one $\psi^R$-factor, which may be
permuted according to the equation ${(\psi^L\psi^R a)}$
strictified.

\section{Coherence with linear endo\-functors}

A symmetric monoidal category is, as usual, a monoidal category
with the natural isomorphism whose components are the arrows
\[
c_{A,B}\!:A\otimes B\str B\otimes A,
\]
which satisfy Mac Lane's coherence conditions (see \cite{ML98},
Section XI.1).

A \emph{linear} endo\-functor in a symmetric monoidal category is
a monoidal endo\-functor $E$ that preserves $c$ globally; i.e.\ we
have the equation
\begin{tabbing}
\hspace{1.5em}\=${(\psi c)}$\hspace{3em}\=$Ec_{A,B}\cirk
\psi_{A,B}=\psi_{B,A}\cirk c_{EA,EB}$.
\end{tabbing}
(We use \emph{linear} instead of \emph{symmetric monoidal} for the
sake of brevity; \emph{linear} comes from the connection with the
structural fragment of \emph{linear} logic, whose name comes from
\emph{linear} algebra.)

A \emph{locally linear} endo\-functor in a symmetric monoidal
category may be defined as a locally monoidal endo\-functor $E$
that satisfies
\begin{tabbing}
\hspace{1.5em}\=${(\psi c)}$\hspace{3em}\=$Ec_{A,B}\cirk
\psi_{A,B}=\psi_{B,A}\cirk c_{EA,EB}$.\kill

\>${(\psi^L\psi^R
c)}$\>$Ec_{A,B}\cirk\psi^L_{A,B}=\psi^R_{B,A}\cirk c_{EA,B}$.
\end{tabbing}
An alternative, simpler, definition is that it is either a left
monoidal or a right monoidal endo\-functor in a symmetric monoidal
category. If it is left monoidal, then from ${(\psi^L\psi^R c)}$
we obtain the definition of $\psi^R$ in terms of $\psi^L$ and $c$,
and we derive ${(\psi^R a)}$, ${(\psi^R l)}$ and ${(\psi^L\psi^R
a)}$.

Let \Mc\ and \Lc\ be the free symmetric monoidal categories with a
family of respectively linear or locally linear endo\-functors;
these categories are defined analogously to $\cal M$ and $\cal L$.
We define the functors $G$ from \Mc\ and \Lc\ to \Fun\ by
stipulating first that $GA$ is the number of occurrences of
generators in $A$. Up to now we took $GA$ to be just the number of
occurrences of generating functors in $A$, but we could as well
have counted also occurrences of generating objects; this was
however superfluous up to now. The remainder of the definitions of
the new functors $G$ is analogous to the definitions of $G$ from
$\cal M$ and $\cal L$ to \Fun, save that we add the clause
corresponding to the picture
\begin{center}
\begin{picture}(30,50)

\put(2,16){\line(4,5){16}}

\put(6,16){\line(4,5){16}}

\put(2,16){\line(1,0){4}}

\put(18,36){\line(1,0){4}}

\put(18,16){\line(-4,5){16}}

\put(22,16){\line(-4,5){16}}

\put(2,36){\line(1,0){4}}

\put(18,16){\line(1,0){4}}

\put(0,42){\makebox(0,0)[l]{$A\otimes B$}}
\put(-1.75,8){\makebox(0,0)[l]{$B\otimes A$}}

\put(-45,25){\makebox(0,0)[l]{$Gc_{A,B}$}}
\end{picture}
\end{center}

The category \Mc\ is equivalent to its strictification \Mcstr, as
$\cal M$ is equivalent to \Mstr\ (see Section~4), and as before we
prove the following proposition, which entails \Mc-Coherence.

\prop{Proposition~2}{The functor $GH$ from \Mcstr\ to \Fun\ is
faithful.}

\dkz We introduce the following abbreviation in \Mcstr:
\[
\Psi_{A_1,A_2;B}=_{df}(\psi_{A_1,A_2}\otimes\mj_B)\cirk(\mj_{EA_1}\otimes
c_{B,EA_2})\!:EA_1\otimes B\otimes EA_2\str E(A_1\otimes
A_2)\otimes B.
\]
We obtain from a developed arrow term of \Mcstr\ a
$\Psi$-\emph{developed} arrow term by replacing the head
$\psi_{A_1,A_2}$ of every $\psi$-factor by $\Psi_{A_1,A_2;I}$;
this replacement is justified by the equation
$\psi_{A_1,A_2}=\Psi_{A_1,A_2;I}$ of \Mcstr. Every
$\Psi$-developed arrow term is equal in \Mcstr\ to an arrow term
of the form $h\cirk g_m\cirk\ldots\cirk g_1$, where in $h$ we have
no occurrences $\Psi$ and $c$, while occurrences of $\cirk$ are
restricted as in the proof of Proposition~1; for
$g_m\cirk\ldots\cirk g_1$ we suppose that it is $\Psi$-developed
without $\psi_0$-factors; i.e.\ it has only $\Psi$-factors and
$c$-factors.

If $m\geq 1$, and $g_i$ is $\Psi$-factor, then we assign to $g_i$
a finite ordinal $\tau(g_i)$ exactly as we did in the proof of
Proposition~1. We then proceed in principle as in that proof to
obtain a normal form. When $\tau(g_{i+1})<\tau(g_i)$ we proceed
exactly as before. Here are the new cases we have to consider.

Suppose we have the $\Psi$-factors $g_i$ and $g_{i+1}$ such that
$\tau(g_{i+1})=\tau(g_i)$. Then we may have the opportunity to
apply the following equations of \Mcstr\ from left to right:
\begin{tabbing}
\hspace{1.5em}\=${(\Psi\Psi 1)}$\hspace{1em}\=$(\Psi_{A_1\otimes
A_3,A_2;B_1}\otimes\mj_{B_2})\cirk\Psi_{A_1,A_3;B_1\otimes
EA_2\otimes B_2}=$
\\*[1ex]
\`$(E(\mj_{A_1}\otimes c_{A_2,A_3})\otimes\mj_{B_1\otimes
B_2})\cirk\Psi_{A_1\otimes A_2,A_3;B_1\otimes
B_2}\cirk(\Psi_{A_1,A_2;B_1}\otimes\mj_{B_2\otimes EA_3})$,
\\[2ex]
\>${(\Psi\Psi 2)}$\>$(\Psi_{A_1,A_2\otimes
A_3;B_1}\otimes\mj_{B_2})\cirk(\mj_{EA_1\otimes
B_1}\otimes\Psi_{A_2,A_3;B_2})=$
\\*[1ex]
\`$\Psi_{A_1\otimes A_2,A_3;B_1\otimes
B_2}\cirk(\Psi_{A_1,A_2;B_1}\otimes\mj_{B_2\otimes EA_3})$.
\end{tabbing}

We call a $c$-factor \emph{atomized} when in its head $c_{A,B}$
the objects $A$ and $B$ are atomic (see Section~2). By the
strictified version of Mac Lane's hexagonal coherence condition
for symmetric monoidal categories (see \cite{ML98}, Section XI.1),
and by $c_{A,I}=l^{-1}_A\cirk r_A=\mj_A$, we may assume that all
our $c$-factors are atomized. Suppose we have an atomic $c$-factor
$g_i$ and a $\Psi$-factor $g_{i+1}$. Then we may have the
opportunity to apply either the naturality and functorial
equations, or the equation ${(\psi c)}$, or the following
equations of \Mcstr:
\begin{tabbing}
\hspace{1.5em}\=${(\Psi c
1)}$\hspace{1em}\=$\Psi_{A_1,A_2;B_1\otimes
B_2}\cirk(c_{B_1,EA_1}\otimes\mj_{B_2\otimes EA_2})=$
\\*[1ex]
\`$(c_{B_1,E(A_1\otimes
A_2)}\otimes\mj_{B_2})\cirk(\mj_{B_1}\otimes\Psi_{A_1,A_2;B_2})$,
\\[2ex]
\>${(\Psi c
2)}$\>$(\mj_{B_1}\otimes\Psi_{A_1,A_2;B_2})\cirk(c_{EA_1,B_1}\otimes\mj_{B_2\otimes
EA_2})=$
\\*[1ex]
\`$(c_{E(A_1\otimes
A_2),B_1}\otimes\mj_{B_2})\cirk\Psi_{A_1,A_2;B_1\otimes B_2}$,
\end{tabbing}
in order to obtain $g_{i+1}\cirk g_i=g'_{i+1}\cirk g'_i$ for a
$\Psi$-factor $g'_i$ and a $c$-factor $g'_{i+1}$. Otherwise, we
must have the opportunity to apply the equations
\begin{tabbing}
\hspace{1.5em}\=${(\Psi c
3)}$\hspace{1em}\=$(\Psi_{A_1,A_2;B_1}\otimes\mj_{B_2})\cirk(\mj_{EA_1\otimes
B_1}\otimes c_{B_2,EA_2})=\Psi_{A_1,A_2;B_1\otimes B_2}$,
\\[2ex]
\>${(\Psi c 4)}$\>$\Psi_{A_1,A_2;B_1\otimes
B_2}\cirk(\mj_{EA_1\otimes B_1}\otimes
c_{EA_2,B_2})=\Psi_{A_1,A_2;B_1}\otimes\mj_{B_2}$,
\end{tabbing}
which follow from the definition of $\Psi$, in order to obtain
$g_{i+1}\cirk g_i=g$ for a $\Psi$-factor $g$. By applying all
these reductions we reach our normal form, which looks as follows.

Let a \emph{block} $\vec{f}_i$ be a composition of $\Psi$-factors
$f_{i_k}\cirk\ldots\cirk f_{i_1}$, all with the same $\tau$ value,
and such that $f_{i_{j+1}}\cirk f_{i_j}$ is never of the form of
the left-hand side of ${(\Psi\Psi 1)}$ and ${(\Psi\Psi 2)}$. Our
normal form is $h\cirk g\cirk\vec{f}_n\cirk\ldots\cirk\vec{f}_1$
such that $n\geq 0$ and the arrangement of the blocks is strictly
increasing (see the proof of Proposition~1 in Section~4); the
arrow term $g$ has no occurrence of $\psi$ and $\psi_0$ (but $c$
may occur), and $h$ has no occurrence of $\psi$ and $c$ (but
$\psi_0$ may occur); $\cirk$ may occur in $h$ only as specified in
the proof of Proposition~1.

The last part of the proof is obtained with slight modifications
of the last part of the proof of Proposition~1. We have the same
kind of induction, but in the basis we do not have just $h=h'$,
but $h\cirk g=h'\cirk g'$. That $h=h'$ follows as before, while
$g=g'$ follows by a coherence result generalizing Mac Lane's
symmetric monoidal coherence (see \cite{ML98}, Section XI.1, or
\cite{DP04}, Chapter~5) as $\cal E$-Coherence of Section~2
generalizes Mac Lane's monoidal coherence. This result is proved
analogously to $\cal E$-Coherence. \qed

\vspace{2ex}

\noindent From this proposition we infer \Mc-Coherence.

To prove \Lc-Coherence we proceed as for \Mc-Coherence. We
introduce the strictification \Lcstr\ of \Lc\ and we prove the
following proposition, from which we will infer \Lc-Coherence.

\prop{Proposition~3}{The functor $GH$ from \Lcstr\ to \Fun\ is
faithful.}

\dkz We have a normal form for the arrow terms of \Lcstr\ which is
a modification of the normal form of the proof of Proposition~2.
For this normal form we have the following abbreviations in
\Lcstr:
\begin{tabbing}
\hspace{1.5em}\=$\Psi^L_{A_1,A_2;B}=_{df}(\psi^L_{A_1,A_2}\otimes\mj_B)\cirk(\mj_{EA_1}\otimes
c_{B,A_2})\!:$\\* \`$EA_1\otimes B\otimes A_2\str E(A_1\otimes
A_2)\otimes B$,
\\[2ex]
\>$\Psi^R_{A_1,A_2;B}=_{df}(\psi^R_{A_1,A_2}\otimes\mj_B)\cirk(\mj_{A_1}\otimes
c_{B,EA_2})\!:$\\* \`$A_1\otimes B\otimes EA_2\str E(A_1\otimes
A_2)\otimes B$,
\end{tabbing}
which are both obtained from the definition of $\Psi_{A_1,A_2;B}$
by adding the superscripts $L$ or $R$ to $\psi$ and deleting some
occurrences of $E$ in the subscripted indices.

The $h$ part of the normal form with $\psi_0$-factors does not
exist now, and instead of $\Psi$-factors we have $\Psi^L$-factors
and $\Psi^R$-factors, which we call collectively $\Psi$-factors.
For a $\Psi$-factor $g_j$ we define $\tau(g_j)$ as before, and we
proceed as before when $\tau(g_{i+1})<\tau(g_i)$ for the
$\Psi$-factors $g_i$ and $g_{i+1}$.

When we have $\tau(g_{i+1})=\tau(g_i)$, then we may apply one of
the following equations of \Lcstr\ from left to right:
\begin{tabbing}
\hspace{.5em}\=${(\Psi^L\Psi^L)}$\hspace{3em}\=$(\Psi^L_{A_1\otimes
A_3,A_2;B_1}\otimes\mj_{B_2})\cirk\Psi^L_{A_1,A_3;B_1\otimes
A_2\otimes B_2}=$
\\*[1ex]
\`$(E(\mj_{A_1}\otimes c_{A_2,A_3})\otimes\mj_{B_1\otimes
B_2})\cirk\Psi^L_{A_1\otimes A_2,A_3;B_1\otimes
B_2}\cirk(\Psi^L_{A_1,A_2;B_1}\otimes\mj_{B_2\otimes A_3})$,
\\[2ex]
\>${(\Psi^L\Psi^R 1)}$\>$\Psi^L_{A_1\otimes A_2,A_3;B_1\otimes
B_2}\cirk(\Psi^R_{A_1,A_2;B_1}\otimes\mj_{B_2\otimes A_3})=$
\\*[1ex]
\`$(\Psi^R_{A_1,A_2\otimes
A_3;B_1}\otimes\mj_{B_2})\cirk(\mj_{A_1\otimes
B_1}\otimes\Psi^L_{A_2,A_3;B_2})$,
\\[2ex]
\>${(\Psi^L\Psi^R 2)}$\>$(\Psi^L_{A_1\otimes
A_3,A_2;B_1}\otimes\mj_{B_2})\cirk\Psi^R_{A_1,A_3;B_1\otimes
A_2\otimes B_2}=$
\\*[1ex]
\`$(E(\mj_{A_1}\otimes c_{A_2,A_3})\otimes\mj_{B_1\otimes
B_2})\cirk(\Psi^R_{A_1,A_2\otimes
A_3;B_1}\otimes\mj_{B_2})\cirk(\mj_{A_1\otimes
B_1}\otimes\Psi^R_{A_2,A_3;B_2})$.
\end{tabbing}

The equation ${(\Psi^L\Psi^L)}$ is obtained from the equation
${(\Psi\Psi 1)}$ in the proof of Proposition~2 above by adding the
superscripts $L$ to $\Psi$ and deleting both occurrences of $E$ in
the subscripted indices. The equation ${(\Psi^L\Psi^R 1)}$ is
obtained in a similar manner from ${(\Psi\Psi 2)}$ read from right
to left. The equation ${(\Psi^L\Psi^R 2)}$, which is analogous to
${(\Psi^L\Psi^L)}$, could be obtained similarly from an equation
of \Mcstr, which we did not need, and did not mention before.

We have moreover eight equations obtained from the equations
${(\Psi c 1)}$-${(\Psi c 4)}$ by adding uniformly the superscripts
$L$ or $R$ to $\Psi$ and deleting some occurrences of $E$ in the
subscripted indices. These equations enable us to obtain a normal
form that looks as follows.

A \emph{block} $\vec{f}_i$ is a composition of $\Psi$-factors
$f_{i_k}\cirk\ldots\cirk f_{i_1}$ all with the same $\tau$ value,
such that $f_{i_{j+1}}\cirk f_{i_j}$ is never of the form of the
left-hand side of ${(\Psi^L\Psi^L)}$, and it is never the case
that $f_{i_j}$ is a $\Psi^R$-factor while $f_{i_{j+1}}$ is a
$\Psi^L$-factor. Our normal form is
$g\cirk\vec{f}_n\cirk\ldots\cirk \vec{f}_1$ such that, as before,
$n\geq 0$ and the arrangement of the blocks is strictly increasing
(see the proof of Proposition~1); the arrow term $g$ has no
occurrence of $\psi^L$ and $\psi^R$ (but $c$ may occur). With this
normal form we proceed as in the proofs of Propositions~1 and 2.
\qed

\vspace{2ex}

\noindent If we define $\Psi^R_{A_1,A_2;B}$ as
\[
(\mj_B\otimes\psi^R_{A_1,A_2})\cirk(c_{A_1,B}\otimes\mj_{EA_2})\!:A_1\otimes
B\otimes EA_2\str B\otimes E(A_1\otimes A_2),
\]
then $\Psi^L$ and $\Psi^R$ would be more symmetric, and we could
use a modification of our normal form that would not favour
pushing $\Psi^L$ to the right as in ${(\Psi^L\Psi^R 1)}$ and
${(\Psi^L\Psi^R 2)}$. In that case, however, our exposition would
be somewhat less economical.

\section{Coherence with conjunctive relevant endo\-func\-tors}

A \emph{conjunctive relevant} category is a symmetric monoidal
category with a diagonal natural transformation, whose components
are the arrows
\[
\Delta_A\!:A\str A\otimes A.
\]
For $\Delta$ we assume the following coherence equations:
\begin{tabbing}
\hspace{1.5em}\=${(\Delta
a)}$\hspace{3em}\=$a_{A,A,A}\cirk(\Delta_A\otimes\mj_A)\cirk\Delta_A=(\mj_A\otimes\Delta_A)\cirk\Delta_A$,
\\[2ex]
\>${(\Delta l)}$\>$l_I\cirk\Delta_I=\mj_I$,
\\[2ex]
\>${(\Delta c)}$\>$c_{A,A}\cirk\Delta_A=\Delta_A$,
\\[2ex]
with $c^m_{A,B,C,D}=_{df} a^{-1}_{A,C,B\otimes
D}\cirk(\mj_A\!\otimes\!(a_{C,B,D}\cirk(c_{B,C}\otimes\mj_D)\cirk
a^{-1}_{B,C,D}))\cirk a_{A,B,C\otimes D}$:\\*[.5ex] \`$(A\otimes
B)\otimes(C\otimes D)\str (A\otimes C)\otimes(B\otimes D)$,
\\*[1ex]
\>${(\Delta ac)}$\>$\Delta_{A\otimes B}=
c^m_{A,A,B,B}\cirk(\Delta_A\otimes\Delta_B)$,
\end{tabbing}
(see \cite{DP96}, Section~2, \cite{P02}, Section~1, and
\cite{DP04}, Sections 9.1-2; the denomination \emph{relevant}
comes from the connection with the structural fragment of
\emph{relevant} logic).

A \emph{conjunctive relevant} endo\-functor in a conjunctive
relevant category is a linear endo\-functor $E$ in this category
that preserves $\Delta$ globally; i.e.\ we have the equation
\begin{tabbing}
\hspace{1.5em}\=${(\Delta a)}$\hspace{3em}\=\kill
\>${(\psi\Delta)}$\>$E\Delta_A=\psi_{A,A}\cirk\Delta_{EA}$.
\end{tabbing}

Let $\cal R$ be the free conjunctive relevant category with a
family of conjunctive relevant endo\-functors. We define the
functor $G$ from $\cal R$ to the category \Rel\ of relations
between finite ordinals as the functor $G$ from \Mc\ to \Fun\ with
an additional clause that corresponds to the following picture:
\begin{center}
\begin{picture}(30,35)(0,7)

\put(4,16){\line(1,3){6.5}}

\put(8,16){\line(1,3){6.5}}

\put(4,16){\line(1,0){4}}

\put(10.5,35.5){\line(1,0){4}}

\put(17,16){\line(-1,3){6.5}}

\put(21,16){\line(-1,3){6.5}}

\put(17,16){\line(1,0){4}}

\put(9,42){\makebox(0,0)[l]{$A$}}
\put(-1.75,8){\makebox(0,0)[l]{$A\otimes A$}}

\put(-45,25){\makebox(0,0)[l]{$G\Delta_A$}}
\end{picture}
\end{center}

Let \Rmin\ be the free conjunctive relevant category, and let $G$
from \Rmin\ to \Rel\ (as a matter of fact, $\mbox{\Fun}^{op}$) be
defined by restricting $G$ from $\cal R$ to \Rel. Then one can
find in \cite{P02} (Section~5) a proof of \Rmin-Coherence.

The category $\cal R$ is equivalent to its strictification \Rstr,
as $\cal M$ is equivalent to \Mstr\ (see Section~4), and, as
before, our goal is to prove the following proposition, which
entails $\cal R$-Coherence.

\prop{Proposition~4}{The functor $GH$ from \Rstr\ to \Rel\ is
faithful.}

We prove first the following auxiliary lemma concerning \Mc\ (see
Section~2 for the notions of diversification and scope).

\prop{\Mc-Theoremhood Lemma}{For $A$ diversified on generating
objects and $B$ diversified, there is an arrow $f\!:A\str B$ of
\Mc\ iff the generators of $A$ and $B$ coincide, and for every
generating functor $E^i$ of $B$ the union of the scopes of the
occurrences of $E^i$ in $A$ is equal to the scope of $E^i$ in
$B$.}

\dkz From left to right the lemma is trivially proved by induction
on the length of $f$ in developed form. For the other direction,
suppose $\{E^1,\ldots,E^n\}$ is the set of generating functors of
$B$. We proceed by induction on $n$. If $n=0$, the set we just
mentioned is empty, and we have trivially an arrow from $A$ to $B$
of symmetric monoidal categories.

For the induction step, let $E^1$ in $E^1B_1$ be the leftmost $E$
of $B$. Since $E^1$ is not in the scope of any other $E$ in $B$,
by the assumptions of the lemma, it is not in the scope of any
other $E$ in $A$ either. So we may assume that $A$ is of the form
\[
D_1\otimes E^1A_1\otimes\ldots\otimes D_n\otimes E^1A_n\otimes
D_{n+1},
\]
with parentheses associated arbitrarily, and $D_i$ being
$E^1$-free. It is clear that we have an arrow of \Mc\ from $A$ to
$E^1(A_1\otimes\ldots\otimes A_n)\otimes D_1\otimes\ldots\otimes
D_{n+1}$. By the induction hypothesis there is an arrow of \Mc\
from $A_1\otimes\ldots\otimes A_n$ to $B_1$, and hence an arrow of
\Mc\ from $E^1(A_1\otimes\ldots\otimes A_n)$ to $E^1B_1$. By
appealing again to the induction hypothesis, we have an arrow of
\Mc\ from $p\otimes D_1\otimes\ldots\otimes D_n$ to $B$ in which
$E^1B_1$ is replaced by $p$. From all that we obtain an arrow of
\Mc\ from $A$ to $B$. \qed

\vspace{2ex}

Let $B$ be a part (proper or not) of an object of \Mcstr, denoted
by $A[B]$, and let $A[B']$ be obtained from this object by
replacing $B$ by $B'$. (We replace a single part $B$ by a single
$B'$.) For $f\!:B\str B'$, let $A[f]\!:A[B]\str A[B']$ be
constructed out of $f$ with identity arrows, $\otimes$ and $E$ in
the obvious way. We can then prove the following.

\prop{Lemma~1}{For the arrow term \[ f\!:A[EA_1\otimes D\otimes
EA_2]\str B\] of \Mcstr\ and $g\!:B\str C$ a $\psi$-factor such
that the ordinals corresponding to the outermost occurrences of
$E$ in $EA_1$ and $EA_2$ are respectively $i$ and $j$, and
$(GHf)(i)\neq(GHf)(j)$, while $(GH(g\cirk f))(i)=(GH(g\cirk
f))(j)$, there exists an arrow term \[ f'\!:A[E(A_1\otimes
A_2)\otimes D]\str C\] of \Mcstr\ such that $g\cirk f=f'\cirk
A[\Psi_{A_1,A_2;D}]$.}

\dkz Note first that every arrow term of \Mcstr\ is a substitution
instance of an arrow term of \Mcstr\ with a diversified target. So
we may assume that $C$ in the lemma is diversified. That $f'$
exists follows from the assumption  that we have $g\cirk f$ and
from the \Mc-Theoremhood Lemma. That $g\cirk f=f'\cirk
A[\Psi_{A_1,A_2;D}]$ follows from \Mc-Coherence. \qed

\prop{Remark}{Consider the arrow term $f\!:A[D\otimes D]\str B$ of
\Mcstr\ with $D$ atomic, and let $i$ and $j$, with $i<j$, be the
ordinals corresponding respectively either to the outermost
occurrences of $E$ in $D$, when $D$ is of the form $ED'$, or
otherwise to the two occurrences of $D$, when $D$ is a
propositional letter. If $(GHf)(j)<(GHf)(i)$, then for $h$ being
$f\cirk A[c_{D,D}]$ we have $(GHh)(i)<(GHh)(j)$ and $f=h\cirk
A[c_{D,D}]$.}

A developed arrow term made only of $\Delta$-factors is called a
$\Delta$-\emph{term}. If $h$ is a $\Delta$-term, then for
$GHh\!:GHA\str GHB$ the converse relation $(GHh)^{-1}\!:GHB\str
GHA$ is an onto function. A $\Delta$-term is \emph{atomized} when
for every $\Delta$-factor in it, in the head $\Delta_A$ of this
$\Delta$-factor, $A$ is atomic.

Let $h$ be an atomized $\Delta$-term such that $E^j$ occurs
exactly once in its source. By naturality  and functorial
equations $h$ is equal to an atomized $\Delta$-term $h_3\cirk
h_2\cirk h_1$ such that, for every factor of $h_1$, its head is
neither in the scope of $E^j$ nor is it of the form
$\Delta_{E^jA'}$; all the heads of the factors of $h_2$ are of the
form $\Delta_{E^jA'}$; and all the heads of the factors of $h_3$
are in the scope of $E^j$ (c.f.\ the proof of $\cal E$-Coherence
in Section~2). Analogously, we can transform every atomized
$\Delta$-term into the \emph{normal form} $h_3\cirk h_2\cirk h_1$
\emph{relative to an occurrence} $E^i$ in its source. (One has to
replace this particular occurrence of $E^i$ by a genuinely new
$E^j$, and then factor the newly obtained arrow term as above; at
the end, one substitutes $E^i$ for $E^j$ everywhere in the term.)

A $\Delta$-\emph{capped} arrow term is an arrow term $f\cirk h$ of
\Rstr\ such that $h\!:D\str A$ is a $\Delta$-term and $f\!:A\str
C$ an arrow term of \Mcstr. A $\Delta$-capped arrow term $f\cirk
h$ is \emph{atomized} when $h$ is an atomized $\Delta$-term.

A \emph{short circuit} in an arrow term $f\cirk h$, with
$h\!:D\str A$, is a pair of ordinals $(i,j)$ such that $i,j\in
GHA$, $i<j$, $(GHh)^{-1}(i)=(GHh)^{-1}(j)$ and
$(GHf)(j)=(GHf)(i)$.

A \emph{useless crossing} is defined analogously to a short
circuit save that we have $(GHf)(j)<(GHf)(i)$. For an example of a
short circuit and a useless crossing see the picture after
Lemma~5.

\prop{Lemma~2}{Every arrow term of \Rstr\ is equal to an atomized
$\Delta$-capped arrow term.}

\noindent To prove this lemma we just apply naturality and
functorial equations together with $\Delta_I=\mj_I$, which is
${(\Delta l)}$ strictified, and ${(\Delta ac)}$.

\prop{Lemma~3}{Every arrow term of \Rstr\ is equal to an atomized
$\Delta$-capped arrow term without short circuits.}

\dkz We apply first Lemma~2, and then we proceed by induction on
the number $n$ of short circuits in an atomized $\Delta$-capped
arrow term. If $n=0$, then we are done. If $n>0$, then our arrow
term is of the form $k\cirk g\cirk f\cirk h$, for $f\cirk h$ an
atomized $\Delta$-capped arrow term without short circuits,
$g\cirk f\cirk h$ an atomized $\Delta$-capped arrow term with a
single short circuit $(i,j)$ and $g$ a $\psi$-factor. Let
$h_3\cirk h_2\cirk h_1$ be the normal form of the $\Delta$-term
$h$ relative to the occurrence $E^i$ in the source of $h$ that
corresponds to the ordinal $(GHh)^{-1}(i)$, which is equal to
$(GHh)^{-1}(j)$, and let $g\cirk f$ be transformed according to
Lemma~1; here $h_2$ is not an identity arrow. Now we can apply
naturality and functorial equations to ``permute'' $h_3$ with
$A[\Psi_{A_1,A_2;D}]$, and then the equations ${(\Delta a)}$
strictified, ${(\Delta c)}$ and ${(\psi\Delta)}$ in order to
decrease $n$. After applying ${(\psi\Delta)}$, we may have to
apply again $\Delta_I=\mj_I$ and ${(\Delta ac)}$ in order to
atomize the resulting $\Delta$-capped arrow term to which we apply
the induction hypothesis. \qed

\prop{Lemma~4}{If for the $\Delta$-terms $f,g\!:A\str B$ we have
$GHf=GHg$, then $f=g$.}

\dkz We proceed essentially as in the proof of $\cal E$-Coherence
in Section~2, by relying on \Rmin-Coherence, instead of Mac Lane's
monoidal coherence, in case $E$ does not occur in $A$. When $E$
occurs in $A$, a difference with the proof of $\cal E$-Coherence
is that the interpolant $C$ is determined not only by $A'$ and
$B'$, but we must take into account $Gf$, which is equal to $Gg$.
From $Gf=Gg$, we may infer also $Gf_1=Gg_1$ and $Gf_2=Gg_2$.

Another difference with the proof of $\cal E$-Coherence is that
the number of occurrences of $E$ in $A$ is not equal to this
number in $B$. In $C$ we may have more than one occurrence of
$E^j$. There will however be no essential difference with the
previous proof, because we do not have to deal with $\Delta$-terms
like $E^j\Delta_{E^jA}$, which have $\Delta$ in the scope of $E^j$
and $E^j$ in the index of $\Delta$. Between two occurrences of
$E^j$ in $C$ there will always be a $\otimes$ in whose scope they
are. Hence in $f_1$ and $g_1$ the subformula $E^jD$ will again
amount to a generating object, and in $f_2$ and $g_2$ only what is
within the scope of $E^j$ counts. Since there may be more than one
occurrence of $E^j$ in $C$ we may need to apply the induction
hypothesis more than once to establish that ${f_2=g_2}$. \qed

\vspace{2ex}

Note that the normal form of $f$ relative to an occurrence of
$E^i$ in $A$ is just a refinement of the factorization $f_2$\cirk
$f_1$ used in the proof of $\cal E$-Coherence in Section~2, and in
the proof we have just finished. For $\cal E$-Coherence, we could
take the factorization $f_2\cirk f_1$ to be such that all the
heads of the factors of $f_1$ are in the scope of $E^j$, and no
head of the factors of $f_2$ is in the scope of $E^j$, while in
the case of $\Delta$-terms, switching the roles of $f_1$ and $f_2$
is not possible.

An arrow term of \Rstr\ is in \emph{normal form} when it is an
atomized $\Delta$-capped arrow term without short circuits and
without useless crossings. We can then prove the following.

\prop{Lemma~5}{Every arrow term of \Rstr\ is equal to an arrow
term in normal form.}

\noindent We just apply Lemma~3, the Remark and the equation
${(\Delta c)}$. Note that we could not apply the Remark without
previously applying Lemma~3. For example, we could have
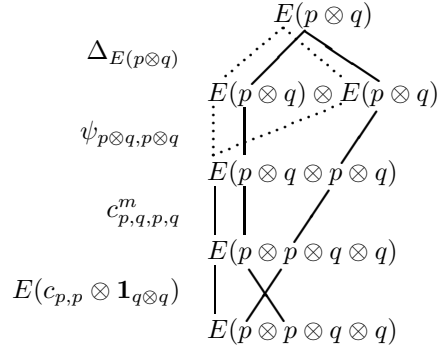
\begin{center}
\begin{picture}(70,140)
\put(3,17){\line(0,1){18}} \put(3,47){\line(0,1){18}}

\multiput(7.5,107)(2.4,2){10}{\makebox(0,0){\circle*{.5}}}
\multiput(51,107)(-2.4,2){10}{\makebox(0,0){\circle*{.5}}}

\multiput(3,77)(0,3){7}{\makebox(0,0){\circle*{.5}}}
\multiput(3,77)(3,1.2){17}{\makebox(0,0){\circle*{.5}}}

\put(25,130){\makebox(0,0)[l]{$E(p\otimes q)$}}
\put(0,100){\makebox(0,0)[l]{$E(p\otimes q)\otimes E(p\otimes
q)$}}

\put(0,70){\makebox(0,0)[l]{$E(p\otimes q\otimes p\otimes q)$}}
\put(0,40){\makebox(0,0)[l]{$E(p\otimes p\otimes q\otimes q)$}}
\put(0,10){\makebox(0,0)[l]{$E(p\otimes p\otimes q\otimes q)$}}

\put(-10,115){\makebox(0,0)[r]{$\Delta_{E(p\otimes q)}$}}
\put(-10,85){\makebox(0,0)[r]{$\psi_{p\otimes q,p\otimes q}$}}
\put(-10,55){\makebox(0,0)[r]{$c^m_{p,q,p,q}$}}
\put(-10,25){\makebox(0,0)[r]{$E(c_{p,p}\otimes\mj_{q\otimes
q})$}}

\thicklines \put(14,47){\line(0,1){18}}
\put(14,77){\line(0,1){18}} \put(15,14){\line(2,3){13}}
\put(29,14){\line(-2,3){13}} \put(33,44){\line(2,3){13}}
\put(51,74){\line(2,3){13}} \put(17,104){\line(1,1){20}}
\put(65,105){\line(-3,2){28}}

\end{picture}
\end{center}
where dotted lines are tied to a short circuit and bold lines to a
useless crossing.

We can now prove Proposition~4.

\vspace{2ex}

\noindent {\sc Proof of Proposition~4. } Note that for an arrow
term $f\cirk h$ of \Rstr\ in normal form, where $h\!:A\str B$ is a
$\Delta$-term and $f\!:B\str C$ is an arrow term of \Mcstr, we
have that $G(f\cirk h)$ determines uniquely $Gh$, $Gf$ and $B$.
This matter, which is not entirely trivial, is established along
the lines of the more general Decomposition Proposition of
\cite{DP09b} (Section~8). We conclude the proof of the proposition
by using Lemma~4 and Proposition~2, i.e.\ \Mc-Coherence. \qed

\section{Coherence in cartesian
categories with relevant endo\-functors }

A \emph{cartesian} category is a conjunctive relevant category
with the monoidal unit object $I$ being a terminal object. The
unique arrow from $A$ to $I$ is $\esp_A\!:A\str I$. This notion of
cartesian category is equivalent to the usual notion, where a
cartesian category is a category with all finite products (see
\cite{DP04}, Sections 9.1-2; some authors use the denomination
\emph{cartesian} for categories with different finite limits than
just finite products; see \cite{Johnst02}, Vol.\ I , Section A1.2,
and \cite{Lei03}, Section 4.1).

In accordance with what we had before, a cartesian endo\-functor
in a cartesian category should preserve $\esp$, which would yield
the equation
\begin{tabbing}
\hspace{1.5em}\=${(\psi\esp)}$\hspace{3em}\=$E\esp_A=\psi_0\cirk\esp_{EA}$,
\end{tabbing}
analogous to ${(\psi\Delta)}$. However, since $GI=\pr$, we see
easily that no definition of $G\esp_A$ would enable us to obtain
coherence with ${(\psi\esp)}$; even the functoriality of $G$,
i.e.\ $G(g\cirk f)=Gg\cirk Gf$, would fail (cf.\ the last picture
in Section~3, and \cite{DP09b}, beginning of Section~5).

We still obtain coherence however for conjunctive relevant
endo\-functors in cartesian categories, and we are going to prove
this now. Let $\cal C$ be the free cartesian category with a
family of conjunctive relevant functors. We define the functor $G$
from $\cal C$ to \Rel\ as the functor $G$ from $\cal R$ to \Rel\
with an additional clause that says that $G\esp$ is the empty
relation between $GA$ and $\pr$, which is $GI$.

The category $\cal C$ is equivalent to its strictification \Cstr,
as $\cal M$ is equivalent to \Mstr\ (see Section~4), and, as
before, our goal is to prove the following proposition, which
entails $\cal C$-Coherence.

\prop{Proposition~5}{The functor $GH$ from \Cstr\ to \Rel\ is
faithful.}

\dkz We proceed analogously to what we had for the proof of
Proposition~4 of the preceding section. We modify the lemmata and
the terminology given there in order to take into account the
presence of $\esp$ in \Cstr.

A developed arrow term made only of $\Delta$-factors and
$\esp$-factors is called a $\Delta\esp$-\emph{term}. If $h$ is a
$\Delta\esp$-term, then for $GHh\!:GHA\str GHB$ the converse
relation $(GHh)^{-1}\!:GHB\str GHA$ is a function. A
$\Delta\esp$-term is \emph{atomized} when for every
$\Delta$-factor and every $\esp$-factor in it, in the heads
$\Delta_A$ or $\esp_A$ of this $\Delta$-factor or $\esp$-factor,
$A$ is atomic.

For every atomized $\Delta\esp$-term $h$, and every occurrence
$E^i$ in its source, we define the \emph{normal form} $h_3\cirk
h_2\cirk h_1$ of $h$ \emph{relative to this occurrence of} $E^i$
exactly as it is defined for atomized $\Delta$-terms in the
preceding section. A $\Delta\esp$-\emph{capped} arrow term is an
arrow term $f\cirk h$ of \Cstr\ such that $h$ is a
$\Delta\esp$-term and $f$ is an arrow term of \Mcstr. A
$\Delta\esp$-capped arrow term $f\cirk h$ is \emph{atomized} when
$h$ is an atomized $\Delta\esp$-term. The notions of short circuit
and useless crossing are defined exactly as in the preceding
section.

An arrow term of \Cstr\ is in \emph{normal form} when it is an
atomized $\Delta\esp$-capped arrow term without short circuits and
without useless crossings. Lemmata 2-5 with \Rstr\ replaced by
\Cstr\ and $\Delta$ by $\Delta\esp$ can be proved along the lines
of the proofs in the preceding section. For the proof of the
modification of Lemma~4, where we relied before on
\Rmin-Coherence, we rely now on cartesian coherence, i.e.\ the
faithfulness of $G$ from the free cartesian category into \Rel\
(see \cite{DP04}, Section 9.2, and references therein).

In contradistinction to what we had in the proof of Proposition~4
at the end of the preceding section, we do not have now any
difficulty in obtaining that $G(f\cirk h)$ determines uniquely
$Gh$, $Gf$ and $B$. This is because we may assume that the target
of $f\cirk h$ is $\otimes$-free. To obtain that, we may compose
with
\begin{tabbing}
\hspace{9em}\=$E^n(r_a\cirk(\mj_a\otimes\esp_B))\;$\=$:E^n(A\otimes
B)\str E^nA$ \quad and\\*
\>$E^n(l_B\cirk(\esp_A\otimes\mj_B))$\>$:E^n(A\otimes B)\str
E^nB$,
\end{tabbing}
for $E^n$ being the sequence of $n$ occurrences of $E$; then we
use the following ``extensionality'' equation of $\cal C$:
\begin{tabbing}
\hspace{1em}$E^{n-1}\psi_{A,B}\cirk\ldots\cirk\psi_{E^{n-1}A,E^{n-1}B}\cirk(E^n(r_a\cirk(\mj_a\otimes\esp_B))
\otimes E^n(l_B\cirk(\esp_A\otimes\mj_B)))\cirk$\\*
\`$\Delta_{E^n(A\otimes B)}=\mj_{E^n(A\otimes B)}$.
\end{tabbing}
Since we may assume that the target of $f\cirk h$ is
$\otimes$-free, we may assume that $h$ is $\Delta$-free, and $Gh$,
$Gf$ and $B$ are then determined uniquely out of $G(f\cirk h)$ in
a straightforward manner. \qed

\section{Coherence in cocartesian categories with endo\-functors}

A \emph{disjunctive relevant} category is a symmetric monoidal
category with a codiagonal natural transformation whose components
are the arrows
\[
\nabla_A\!:A\otimes A\str A.
\]
For $\nabla$ we assume coherence equations dual to ${(\Delta a)}$,
${(\Delta l)}$, ${(\Delta c)}$ and ${(\Delta ac)}$.

A \emph{disjunctive relevant} endo\-functor in a disjunctive
relevant category is a linear endo\-functor $E$ in this category
that satisfies the equation
\begin{tabbing}
\hspace{1.5em}\=${(\psi\esp)}$\hspace{3em}\=\kill
\>${(\psi\nabla)}$\>$E\nabla_A\cirk\psi_{A,A}=\nabla_{EA}$.
\end{tabbing}
This equation, together with others, enables us to reduce to a
propositional letter all the indices of $\nabla$, and that
together with the naturality of $\nabla$ is all we need
essentially to push every occurrence of $\nabla$ to the left. This
enables us to prove coherence for disjunctive relevant categories.

A \emph{cocartesian} category is a disjunctive relevant category
with the monoidal unit object $I$ being an initial object. The
unique arrow from $I$ to $A$ is $!_A\!:I\str A$. Equivalently,
cocartesian categories are defined as categories with all finite
coproducts.

It will follow from the coherence result below that any
endo\-functor $E$ in a cocartesian category is a disjunctive
relevant endo\-functor with the definitions
\begin{tabbing}
\hspace{3.5em}\=$\psi_{A,B}=_{df}\nabla_{E(A\otimes
B)}\cirk(E((\mj_A\otimes\;!_B)\cirk r^{-1}_A)\otimes
E((!_A\otimes\mj_B)\cirk l^{-1}_B))$,\\[1ex]
\>$\psi_0=_{df}\; !_{EI}$.
\end{tabbing}
Moreover, this endo\-functor preserves $!$, in the sense that it
satisfies the equation
\begin{tabbing}
\hspace{1.5em}\=${(\psi\esp)}$\hspace{3em}\=\kill
\>${(\psi!)}$\>$E!_A\cirk\psi_0=\;!_{EA}$.
\end{tabbing}
Unlike ${(\psi\esp)}$, this equation is in accordance with
coherence. Note that the equations of the definitions of $\psi_0$
and $\psi$ above follow from the initiality of $I$, from the
requirement that $\otimes$ is a coproduct and from the equations
${(\psi l)}$ and ${(\psi r)}$.

Let $\cal D$ be the free cocartesian category with a family of
endo\-functors. We define the functor $G$ from $\cal D$ to \Fun\
as the functor $G$ from \Mc\ to \Fun\ with the additional clause
that corresponds to the following picture:
\begin{center}
\begin{picture}(30,35)(0,7)

\put(4,35.5){\line(1,-3){6.5}}

\put(8,35.5){\line(1,-3){6.5}}

\put(4,35.5){\line(1,0){4}}

\put(10.5,16){\line(1,0){4}}

\put(17,35.5){\line(-1,-3){6.5}}

\put(21,35.5){\line(-1,-3){6.5}}

\put(17,35.5){\line(1,0){4}}

\put(9,8){\makebox(0,0)[l]{$A$}}
\put(0,42){\makebox(0,0)[l]{$A\otimes A$}}

\put(-45,25){\makebox(0,0)[l]{$G\nabla_A$}}
\end{picture}
\end{center}
and the clause that says that $G!_A$ is the empty function from
$\pr$, which is $GI$, to $GA$. Then we can prove $\cal
D$-\emph{Coherence}.

We proceed essentially as in the proof of $\cal E$-Coherence in
Section~2, and as in the proofs of Lemma~4 and its modifications
in the two preceding sections. We rely now on cocartesian
coherence, instead of monoidal coherence, \Rmin-Coherence and
cartesian coherence respectively. (The only difference is that now
we work relative to an occurrence $E^j$ in the target $B$, and the
factorization $f_2\cirk f_1$ is such that all the heads of the
factors of $f_1$ are in the scope of $E^j$ and no head of a factor
of $f_2$ is in the scope of $E^j$.)

An alternative way to prove $\cal D$-Coherence is to rely on the
factorization $f_2\cirk f_1$ such that $f_1$ is an arrow term of
\Lc, and in the developed arrow term $f_2$ every factor is either
a $\nabla$-factor or a $!$-factor such that the index of its head
is a propositional letter. For that we use the equations
${(\psi\nabla)}$ and ${(\psi!)}$ above.

Cocartesian coherence, i.e.\ the faithfulness of $G$ from the free
cocartesian category into \Fun, follows from cartesian coherence
(see \cite{DP04}, Section 9.2, and references therein).
Cocartesian coherence may be proved by relying on a normal form
inspired by Gentzen's cut elimination (see \cite{DP04}, Sections
9.1-2), but we could rely alternatively on a developed strictified
normal form $f_3\cirk f_2\cirk f_1$, where $f_1$ has atomized
$c$-factors only, $f_2$ has atomized $\nabla$-factors only, and
$f_3$ has atomized $!$-factors only.

The coherence results of this paper yield coherence results for
categories with arrows oriented in the opposite direction. In
these categories we do not have monoidal functors with $\psi$ and
$\psi_0$, but \emph{comonoidal} functors with arrows oriented
oppositely to $\psi$ and $\psi_0$.

\vspace{2ex}

\noindent {\small {\it Acknowledgement.} We would like to thank
very much an anonymous referee for a careful reading of our text
and for useful suggestions. Work on this paper was supported by
the Ministry of Science of Serbia (Grants 144013 and 144029).}

\end{document}